\documentclass[12pt]{amsart}
\usepackage{amsxtra,latexsym,amssymb}
\usepackage{epsfig,rotate,amsthm}
\usepackage{hyperref}

\def\qed{{\hfill $\Box$}}

\def\C{{\mathbb C}}

\def\K{{\mathbb K}}

\def\A{\mathcal{A}(\alpha, \beta, \phi)}
\theoremstyle{theorem}
\newtheorem{thm}{Theorem}[section]
\newtheorem{cor}{Corollary}[section]
\newtheorem{prop}{Proposition}[section]

\newtheorem{lem}{Lemma}[section]
\theoremstyle{definition}
\newtheorem{defn}{Definition}[section]
\theoremstyle{remark}

\newtheorem{rem}{\bf Remark}[section]

\begin{document}
\title[Down-up algebras over Polynomial rings]{Down-up algebras defined over a polynomial base ring $\K[t_{1}, \cdots, t_{n}]$}
\author[X. Tang]{Xin Tang}
\address{Department of Mathematics \& Computer Science\\
Fayetteville State University\\
1200 Murchison Road, Fayetteville, NC 28301}
\email{xtang@uncfsu.edu} 
\date{\today}
\begin{abstract}
In this paper, we study a class of down-up algebras $\A$ defined over a polynomial base ring $\K[t_{1}, \cdots, t_{n}]$ and establish several analogous results. We first construct a $\K-$basis for the algebra $\A$. As a result, we prove that the Gelfand-Kirillov dimension of $\A$ is $n+3$ and completely determine the center of $\A$ when $char\K=0$. Then, we prove that the algebra $\A$ is a noetherian domain if and only if $\beta\neq 0$; and $\A$ is Auslander-regular when $\beta \neq 0$. We also prove that the global dimension of $\A$ is $n+3$; and the algebra $\A$ is a prime ring except $\alpha=\beta=\phi=0$. Moreover, we obtain some results on the Krull dimension, isomorphisms, and automorphisms of the algebra $\A$.
\end{abstract}
\maketitle 
\section*{Introduction}
The notion of down-up algebras was introduced by Benkart and Roby in \cite{B,BR}. Let $\mathbb{K}$ be a field and $\alpha, \beta, \gamma \in \mathbb{K}$. The down-up algebra $A(\alpha, \beta, \gamma)$ is defined to be a $\mathbb{K}-$algebra generated by $d, u$ subject to the following two relations:
\begin{eqnarray*}
d^{2}u=\alpha dud +\beta ud^{2}+\gamma d,\\
du^{2}=\alpha udu+\beta u^{2}d+\gamma u.\\
\end{eqnarray*}
The family of down-up algebras $\{A(\alpha, \beta, \gamma)\mid \alpha, \beta, \gamma\in \K\}$ contains the enveloping algebras of both the Lie algebra $sl_{2}$ and the three dimensional Heisenberg Lie algebra $\mathcal{H}$ and their various quantizations. 

Down-up algebras have been extensively studied in the literature and a lot of important properties have been established. The representation theory of down-up algebras was thoroughly investigated in \cite{BR, BW, CM, J, Kul} and the references therein; the centers of down-up algebras over $\C$ were determined in \cite{Kul, Z}. The Hopf algebra structures of down-up algebras were investigated in \cite{BW, KM}; the primitivity of down-up algebras was studied in \cite{J, KK1}; and the primitive ideals of down-up algebras over $\C$ were explicitly constructed in \cite{P, PM}. In particular, it was proved in \cite{KMP} that a down-up algebra $A(\alpha, \beta, \gamma)$ is a noetherian domain if and only if the parameter $\beta \neq 0$; and a noetherian down-up algebra is further Auslander-regular of global dimension $3$. Moreover, it was proved in \cite{KK2} that a non-noetherian down-up algebra also has a global dimension of $3$; and all down-up algebras are prime rings except in the case when the parameters $\alpha, \beta, \gamma$ are all zero. The Krull dimensions of down-up algebras were studied in \cite{BL, BV}. In addition, the representations and centers of down-up algebras defined over fields of positive characteristics were studied in \cite{H1, H2}.

As generalizations of down-up algebras, a class of generalized down-up algebras was defined and studied in \cite{CS}, and some down-up algebras from trees were studied in \cite{KK3}. Recently, the automorphism groups of certain generalized down-up algebras (including some families of down-up algebras) were determined \cite{CL} under restrictions on the parameters; and the factorial classification of generalized down-up algebras was obtained in \cite{LL}. Besides, the study of isomorphism problem for down-up algebras was initiated in \cite{BR} and was substantially investigated in \cite{CM}. Motivated by a question posed in \cite{BR}, the homogenized down-up algebras were also investigated in \cite{C}. Note that the homogenized down-up algebras provide many examples of Artin-Schelter regular algebras of global dimension $4$. The invariant theory of finite group action on down-up algebras has recently been investigated in \cite{KK4, KKZ}.

In this paper, we are going to study a family of down-up algebras $\A$, which is instead defined over a polynomial base ring $\K[t_{1}, \cdots, t_{n}]$. These algebras $\A$ are indeed generalizations of the down-up algebras $A(\alpha, \beta, \gamma)$ defined over a field $\K$. As we shall see, these algebras include the homogenized down-up algebras and the tensor products of down-up algebras $A(\alpha, \beta, \gamma)$ with a polynomial ring $\K[t_{1}, \cdots, t_{n}]$ as typical examples. 

Let $\alpha, \beta\in \K$ be scalars, and let $\phi=\phi(t_{1}, \cdots, t_{n}) \in \K[t_{1}, \cdots, t_{n}]$ be a polynomial in $t_{1}, \cdots, t_{n}$. We define a down-up algebra $\A$ as the $\K-$algebra generated by $u, d, t_{1}, \cdots, t_{n}$ subject to the following relations:
\begin{eqnarray*}
t_{i}t_{j}=t_{j}t_{i}\, \text{for}\, i, j =1, \cdots, n;\\
t_{i}u=ut_{i}, t_{i}d=dt_{i}\,\text{for}\, i, j=1, \cdots, n;\\
d^{2}u=\alpha dud+\beta ud^{2}+\phi (t_{1}, \cdots, t_{n}) d;\\
du^{2}=\alpha udu+\beta u^{2}d+u\phi (t_{1}, \cdots, t_{n}).
\end{eqnarray*}
Obviously, the algebras $\A$ can be viewed as central extensions of the down-up algebras $A(\alpha, \beta, \gamma)$. If the parameter $\phi=\phi (t_{1}, \cdots, t_{n})$ is chosen to be a constant, denoted by $\gamma$, then the algebra $\A$ is indeed the tensor product of the down-up algebra $A(\alpha, \beta, \gamma)$ with the polynomial algebra $\K[t_{1}, \cdots, t_{n}]$. 

The algebras $\A$ share many similar properties with the down-up algebras $A(\alpha, \beta, \gamma)$. Using some classical methods on filtered rings \cite{MR}, we will be able to establish many analogous results for the algebras $\A$ based upon the ones for the down-up algebras $A(\alpha, \beta, \gamma)$. In particular, we will construct a $\K-$basis for the algebra $\A$; as an application, we will determine the center of $\A$ in the case when $char \K=0$, and prove that the Gelfand-Kirillov dimension of $\A$ is $n+3$. Once again, we will prove that the algebra $\A$ is a noetherian domain if and only if $\beta \neq 0$, and $\A$ is Auslander-regular when $\beta\neq 0$, and the global dimension of $\A$ is always $n+3$, and $\A$ is a prime ring except $\alpha=\beta=\phi=0$. Besides, we will obtain some results on the Krull dimension, isomorphisms and automorphisms for the algebra $\A$. 

It is no surprise that we are going to closely follow the approaches used in the investigation of down-up algebras $A(\alpha, \beta, \gamma)$. In the sequel, in order to distinguish the algebras $\A$ from $A(\alpha, \beta, \gamma)$, we will refer to the down-up algebras $A(\alpha, \beta, \gamma)$ as ``classical'' down-up algebras. 

The paper is organized as follows. In Section $1$, we will define the down-up algebra $\A$ and establish some basic properties for $\A$. In Section $2$, we will study the noetherianity, global dimension, and Krull dimension of $\A$. In Section $3$, we will study the automorphisms and isomorphisms of $\A$.

\section{The Down-up algebra $\A$}
In this section, we first state the definition of a down-up algebra $\A$ defined over a polynomial ring $\K[t_{1}, \cdots, t_{n}]$. Then we will establish a variety of basic properties for the algebra $\A$.

Assume that $\K$ is an arbitrary field. Let $\alpha, \beta\in \K$ be scalars and $\phi=\phi (t_{1}, \cdots, t_{n}) \in \K[t_{1}, \cdots, t_{n}]$ be a polynomial in $t_{1}, \cdots, t_{n}$. We have the following definition.
\begin{defn}
The (generalized) down-up algebra $\A$ is defined to be a $\K-$algebra generated by $u, d, t_{1}, \cdots, t_{n}$ subject to the following relations:
\begin{eqnarray*}
t_{i}t_{j}=t_{j}t_{i}\,\text{for}\, i, j =1, \cdots, n;\\
t_{i}u=ut_{i}, \, t_{i}d=dt_{i}\, \text{for}\, i, j=1, \cdots, n;\\
d^{2}u=\alpha dud+\beta ud^{2}+\phi (t_{1}, \cdots, t_{n}) d;\\
du^{2}=\alpha udu+\beta u^{2}d+u\phi (t_{1}, \cdots, t_{n}).
\end{eqnarray*}
\end{defn}

The down-up algebra $\A$ is abstractly defined via generators and relations over the base field $\K$. Thus it is a priori to construct a $\K-$basis for the algebra $\A$. Similar to the result for the ``classical" down-up algebras $A(\alpha, \beta, \gamma)$ in \cite{BR}, we have the following construction of a $\K-$basis for $\A$.
\begin{thm}
The set $\{u^{i}(du)^{j}d^{k}t_{1}^{m_{1}}\cdots t_{n}^{m_{n}}\mid i, j, k, m_{1},\cdots, m_{n}\geq 0\}$ is a $\K-$basis of the algebra $\A$. Furthermore, the algebra $\A$ has a Gelfand-Kirillov dimension of $n+3$.
\end{thm}
{\bf Proof:} We will follow the approach used in \cite{BR}. First of all, it is easy to verify that the algebra $\A$ is spanned by the set $ \{u^{i}(du)^{j}d^{k}t_{1}^{m_{1}}\cdots t_{n}^{m_{n}}\mid i, j, k, m_{1}, \cdots, m_{n}\geq 0\}$ as a $\K-$vector space. Now we show 
that the set $\{u^{i}(du)^{j}t_{1}^{m_{1}}\cdots t_{n}^{m_{n}}\mid i,j, k, m_{1}, \cdots, m_{n}\geq 0\}$ is indeed $\K-$linearly independent. 
To this end, we need to use the Diamond Lemma \cite{Ber}. Namely, in the free algebra $\K\langle u, d, t_{1}, \cdots, t_{n}\rangle$ generated by $u, d, t_{1}, \cdots, t_{n}$, we will assign the degrees to the generators as follows:
\[
deg(u)=deg(d)=\text{max}\{deg(\phi), 1\}, \quad deg(t_{1})=\cdots=deg(t_{n}) = 1.
\]
In addition, we shall order the monomials in $u, du, d, t_{1}, \cdots, t_{n}$ first by their total degrees and then lexicographically pursuant to the prescribed order $u<d<t_{1}<\cdots<t_{n}$. It is obvious to see that there is only one possible ambiguity, namely $(d^{2}u)u=d^{2}u^{2}=(d^{2}u)u$, to be resolved. Using the defining relations of the algebra $\A$, one can easily verify the following: 
\[
(d^{2}u-\alpha dud -\beta ud^{2}-\phi d)u-d(du^{2}-\alpha udu -\beta u^{2}d-\phi u)=\beta (ud^{2}u-du^{2}d).
\]
Since $ud^{2}u< d^{2}u^{2}$ and $du^{2}d< d^{2}u^{2}$, the ambiguity is therefore resolvable. Hence, the set $\{u^{i}(du)^{j}d^{k}t_{1}^{m_{1}}\cdots t_{n}^{m_{n}}\mid i, j, k, m_{1},\cdots, m_{n}\geq 0\}$ consisting of irreducible monomials is a $\K-$basis of the algebra $\A$. Using the basis $\{u^{i}(du)^{j}d^{k}t_{1}^{m_{1}}\cdots t_{n}^{m_{n}}\mid i, j, k, m_{1},\cdots, m_{n}\geq 0\}$, one can use a similar argument as in \cite{BR} to prove that the Gelfand-Kirillov dimension of the algebra $\A$ is indeed $n+3$. 
\qed

As a result of the previous theorem, we have the following corollary.
\begin{cor}
\begin{enumerate}
\item The subalgebra $\K[t_{1}, \cdots, t_{n}]$ of $\A$ generated by $t_{1}, \cdots, t_{n}$ is a polynomial ring; 
\item Any non-zero polynomial $f(t_{1}, \cdots, t_{n})\in \K[t_{1}, \cdots, t_{n}]$ is a regular central element of $\A$.
\end{enumerate}
\end{cor}\qed

As generalizations of the ``classical" down-up algebras $A(\alpha, \beta, \gamma)$, the down-up algebras $\A$ are closely related to $A(\alpha, \beta, \gamma)$ in many perspectives. As a matter of fact, we have the following proposition.
\begin{prop} 
\begin{enumerate} 

\item Let $\mathbb{S}=\{f(t_{1}, \cdots, t_{n})\mid 0\neq f(t_{1}, \cdots, t_{n}) \in \K[t_{1}, \cdots, t_{n}]\}$. Then $\mathbb{S}$ is a multiplicative subset of $\A$; and the localization $A_{\mathbb{S}}$ is isomorphic to a ``classical" down-up algebra $A(\alpha, \beta, \gamma)$ defined over the filed $\K(t_{1}, \cdots, t_{n})$ for some $\gamma \in \K(t_{1}, \cdots, t_{n})$.
\item Let $I= (t_{1}-\lambda_{1}, \cdots, t_{n}-\lambda_{n})$, where $\lambda_{1}, \cdots, \lambda_{n} \in \K$ be an ideal of $\A$. Then the quotient algebra $\A \slash I$ is isomorphic to a ``classical" down-up algebra $A(\alpha, \beta, \gamma)$ defined over the field $\K$ for some $\gamma \in \K$.
\end{enumerate}
\end{prop} 
\qed

When $\K=\C$, the centers of the down-up algebras $A(\alpha, \beta, \gamma)$ were completely determined in \cite{Z} and independently in \cite{Kul} when $\alpha^{2}+4\beta \neq 0$. The centers of the homogenized down-up algebras were also described in \cite{Z}. It is straightforward to verify the method used in \cite{Z} can be applied to a down-up algebra $A(\alpha, \beta, \gamma)$ defined over any sufficiently large field $\K$ of characteristic zero. In particular, we can apply the method to study the center of the algebra $\A_{\mathbb{S}}$. As a result, we have the following description of the center of the down-up algebra $\A$ via pulling back information from $\A_{\mathbb{S}}$.
\begin{thm}
Assume that $char \K=0$ and $\K$ is sufficiently large so that the quadratic equation $x^{2}-\alpha x-\beta=0$ has two roots $r, s$ in $\K$. Then the center $Z(\A)$ of the algebra $\A$ is described as follows.\\

\begin{enumerate}
\item If $r=s$ is a primitive root of unity with order $m\geq 2$, then $Z(\A)=\K[(du-rdu+\frac{\phi}{r-1})^{m}, t_{1}, \cdots, t_{n}]$. If $\phi=0$, then we will allow $m=1$.\\
\item $Z(\mathbb{A}(2, -1, \phi)=\K[(du-ud)^{2}-\phi (du+ud), t_{1}, \cdots, t_{n}]$ if $\phi \neq 0$.\\
\item $Z(\A)=\K[(du-sud+\frac{\phi}{r-1})^{m}, t_{1}, \cdots, t_{n}]$ if $r$ is a primitive root of unity of order $m>1$ and $s$ is not a root of unity. We will allow $m=1$ when $\phi=0$.\\
\item $Z(\A)=\K[(du-rud+\frac{\phi}{s-1})^{m}, t_{1}, \cdots, t_{n}]$ if $s$ is a primitive root of unity of order $m>1$ and $r$ is not a root of unity. We will allow $m=1$ when $\phi=0$.\\
\item $Z(\A)=\K[(du-rud+\frac{\phi}{s-1})^{j}(du-sud+\frac{\phi}{r-1})^{i}\, \text{where}\, r^{i}s^{j}=1, t_{1}, \cdots, t_{n}]$ if $r\neq s$ and neither $r$ nor $s$ is a root of unity.\\

\item $Z(\A)=\K[d^{m}, u^{m}, t_{1}, \cdots, t_{n}, (du-rud+\frac{\phi}{s-1})^{i}(du-sud+\frac{\phi}{r-1})^{j}\, \text{where}\, s^{i}r^{j}=1]$ if $r\neq s$ and $r,s$ are both primtive roots of unity of orders $m_{1}>1, m_{2}>1$ and $m$ is the least common multiple of $m_{1}, m_{2}$. We can allow $m_{1}=1$ or $m_{2}=1$ when $\phi=0$.\\

\item If $\phi \neq 0$ and $r=1$ and $s$ is a primitive root of unity of order $m>1$, then $Z(\A)=\K[(du-rud+\frac{\phi}{s-1})^{m}, t_{1}, \cdots, t_{n}]$.\\
\item If $\phi \neq 0$ and $s=1$ and $r$ is a primitive root of unity of order $m>1$, then $Z(\A)=\K[(du-sud+\frac{\phi}{r-1})^{m}, t_{1}, \cdots, t_{n}]$.\\
\item Otherwise, $Z(\A)=\K[t_{1}, \cdots, t_{n}]$.
\end{enumerate}
\end{thm}

{\bf Proof:} Note $\mathbb{S}=\{f(t_{1}, \cdots, t_{n})\mid 0\neq f(t_{1}, \cdots, t_{n})\in \K[t_{1}, \cdots, t_{n}]\}$ is a multiplicative subset of $\A$ consisting of regular central elements, and the localization $\A_{\mathbb{S}}$ of $\A$ with respect to $\mathbb{S}$ is isomorphic to a ``classical" down-up algebra $A(\alpha, \beta, \gamma)$ defined over the field $\K(t_{1}, \cdots, t_{1})$ for some $\alpha, \beta \in \K, \gamma \in \K(t_{1}, \cdots, t_{n})$. Since $char \K=0$ and the equation $x^{2}-\alpha x-\beta=0$ has two roots in $\K$, we have that $char \K(t_{1}, \cdots, t_{n})=0$ and the equation $x^{2}-\alpha x-\beta=0$ also has two roots in $\K(t_{1}, \cdots, t_{n})$. Now it is straightforward to verify that the method used in \cite{Z} can be adopted to completely determine the center of the ``classical" down-up algebra $\A_{\mathbb{S}}$. As a result, the description on the center of a down-up algebra $A(\alpha, \beta, \gamma)$ over $\C$ as obtained {\bf Theorem 1.3} in \cite{Z} can be directly applied to the algebra $\A_{\mathbb{S}}$. Therefore, we can have the desired result on the center of the algebra $\A$ via pulling the information back to the algebra $\A$ from its localization $\A_{\mathbb{S}}$. \qed

\begin{rem}
When the field $\K$ is of positive characteristic, one can similarly describe the center of $\A$ via pulling back results on the centers of ``classical" down-up algebras over fields of prime characteristics as established in \cite{H2} through the localization. However, we will not repeat the details here.
\end{rem}
\qed

Recall that many classes of (no necessarily irreducible) representations have been constructed for the ``classical" algebras $A(\alpha, \beta, \gamma)$ in the references \cite{BR, BW, J, Kul} and therein; and the primitive ideals of $A(\alpha, \beta, \gamma)$ defined over $\C$ have also been constructed in \cite{P, PM}. We now state a few basic results on the irreducible representations and primtive ideals of the down-up algebras $\A$ in terms of the results on the irreducible representations and primitive ideals of the ``classical" down-up algebras $A(\alpha, \beta, \gamma)$.
\begin{prop}
Assume that $\K$ is an algebraically closed uncountable field. Each irreducible representation $M$ of the algebra $\A$ can be viewed as an irreducible representation of a quotient algebra $\A/I$ (which is indeed a ``classical" down-up algebra defined over $\K$) for some ideal $I=(t_{1}-\lambda_{1}, \cdots, t_{n}-\lambda_{n})$; and any primitive ideal $P$ of the algebra $\A$ is the inverse image of a primitive ideal $P^{\prime}$ of a certain quotient algebra $\A/I$ for some $I=(t_{1}-\lambda_{1}, \cdots, t_{n}-\lambda_{n})$ via the algebra homomorphism $q_{I}\colon \A\longrightarrow A/I$.
\end{prop}
{\bf Proof:} First of all, the algebra $\A$ has a countable $\K-$basis. Then every irreducible representation $M$ of $\A$ has a countable dimension over the base field $\K$. Since the base field $\K$ is an algebraically closed uncountable field, the result now follows from Schur's lemma (see \cite{Dix}). Second of all, note that a primitive ideal $P$ of the algebra $\A$ is the annihilator of an irreducible representation $M$ of $\A$. Since the irreducible representation $M$ is also regarded as an irreducible representation of a certain quotient algebra $\A/I$ of $\A$ for some $I=(t_{1}-\lambda_{1}, \cdots, t_{n}-\lambda_{n})$, we have that $I\subseteq P$ and $P^{\prime}=P/I$ is a primitive ideal of the quotient algebra $\A/I$. Thus we have the desired description of the primitive ideals of $\A$. \qed

\section{Global and Krull dimensions of $\A$}
The ``classical" down-up algebra $A(\alpha, \beta, \gamma)$ enjoys many nice homological properties \cite{KK1, KMP}. Indeed, the ``classical" down-up algebra $A(\alpha, \beta, \gamma)$ is a noetherian domain if and only if $\beta\neq 0$; and when $\beta \neq 0$, the algebra $A(\alpha, \beta, \gamma)$ is also Auslander-regular of global dimension $3$ (see \cite{KMP}). Later on, it was further proved in \cite{KK1} that any non-noetherian ``classical" down-up algebra $A(\alpha, \beta, \gamma)$ also has a global dimension of $3$. The Krull dimensions of noetherian ``classical" down-up algebras were studied in \cite{BL,BV}. In particular, the Krull dimensions of noetherian down-up algebras $A(\alpha, \beta, \gamma)$ were explicitly described in \cite{BL}.

In this section, we will first establish a necessary-sufficient condition for the down-up algebra $\A$ to be a noetherian domain. Then we will calculate the global dimension of any down-up algebra $\A$, and establish some results on the Krull dimension of the algebra $\A$ under some conditions on the parameters. Note that the algebra $\A$ is obviously isomorphic to its opposite ring $\A^{op}$ via the mapping $u\rightarrow d^{op}, d\rightarrow u^{op}, t_{1}\rightarrow t^{op}_{1},\cdots, t_{n}\rightarrow t_{n}^{op}$. Therefore, there is no need to distinguish the left and right global dimensions (resp. Krull dimensions).

First of all, we have the following exact analogue of the result on noetherianity for the ``classical" down-up algebras $A(\alpha, \beta, \gamma)$ as established in \cite{KMP}.
\begin{thm}
The following statements are equivalent.

\begin{enumerate}
\item The algebra $\A$ is left (and right) noetherian;\\

\item The algebra $\A$ is a domain;\\

\item $\beta \neq 0$;\\

\item The subalgebra $\K[ud, du, t_{1}, \cdots, t_{n}]$ is a polynomial ring in $n+2$ indeterminates: $ud, du, t_{1}, \cdots, t_{n}$.
\end{enumerate}
In particular, when $\beta\neq 0$, the algebra $\A$ is also Auslander-regular.
\end{thm}

{\bf Proof:} Suppose $\beta\neq 0$. Following the idea in \cite{KMP}, we can embed the algebra $\A$ into a skew group ring. We will spell out the details. Let $R=\K[x, y, t_{1}, \cdots, t_{n}]$ be a polynomial ring in the indeterminates $x, y, t_{1}, \cdots, t_{n}$ and define an automorphism $\sigma$ of $R$ as follows:
\[
\sigma(x)=y,\quad \sigma(y)=\alpha x+\beta y+\phi,\quad \sigma(t_{i})=t_{i}\quad \text{for}\quad i=1, \cdots, n.
\]
It is obvious that the mapping $\sigma$ is an algebra automorphism of $R$ because $\beta \neq 0$ and $\phi$ is a polynomial in $t_{1}, \cdots, t_{n}$. Let $S=R[z^{\pm 1};\sigma]$ be the corresponding skew group ring of the infinite cyclic group $\langle z\rangle$ over the base ring $R$. Indeed, we will require that $rz=z\sigma(r)$ for any $r\in R$. Let us set $D=z^{-1}$ and $U=xz$ in the ring $S$. Then it is easy to verify that $UD=xzz^{-1}=x$ and $DU=z^{-1}xz=\sigma(x)=y$. As a result, we have $D^{2}U=\alpha DUD+\beta UD^{2}+\phi D$ and $DU^{2}=\alpha UDU +\beta U^{2}D+\phi U$ in $S$. Thus there is an algebra homomorphism $\theta \colon \A\longrightarrow S$ defined by $\theta(d)=D, \theta(u)=U$ and $\theta(t_{i})=t_{i}$ for $i=1, \cdots, n$. Therefore, the elements $ud, du, t_{1}, \cdots, t_{n}$ are algebraically independent over $\K$; and in particular, the subalgebra $\K[ud, du, t_{1}, \cdots, t_{n}]$ is a polynomial ring in $ud, du, t_{1}, \cdots, t_{n}$. One can verify the set $\mathcal{B}=\{(ud)^{i}(du)^{j}d^{k}t_{1}^{m_{1}}\cdots t_{n}^{m_{n}}, (ud)^{i}(du)^{j}du^{k+1}t_{1}^{m_{1}}\cdots t_{n}^{m_{n}}\mid i, j, k, m_{1}, \cdots, m_{n}\geq 0\}$ spans the algebra $\A$ and the image of $\mathcal{B}$ under $\theta$ is the following set:
\[
\{x^{i}y^{j}z^{-k}t_{1}^{m_{1}}\cdots t_{n}^{m_{n}}, x^{i}y^{j}(xz)^{k+1}t_{1}^{m_{1}}\cdots t_{n}^{m_{n}}\mid i, j, k, m_{1}, \cdots, m_{n}\geq 0\}.
\]
Thus we know that $\mathcal{B}$ is indeed a $\K-$basis of $\A$ and $\theta$ is also a monomorphism. Therefore, we have established the implications: $(3)\Rightarrow (2)$, $(3)\Rightarrow (4)$.

Let us maintain the assumptions on $\beta, R, \sigma$ and denote by $R(\sigma, x)$ the generalized Weyl algebra over $R$ generated by the variables $X^{\pm}$ subject to the relations:
\[
X^{-}X^{+}=x,\quad X^{+}X^{-}=\sigma(x),\quad X^{+}r=\sigma(r)X^{+}, \quad X^{-}r=\sigma(r)X^{-}
\]
for all $r\in R$. It is straightforward to verify that $\A$ is isomorphic to the generalized Weyl algebra $R(\sigma, x)$. Thus, $\A$ is both a left and right noetherian ring. So we have proved the implication: $(3)\Rightarrow (1)$.

Conversely, assume that $\beta=0$, then $d(du-\alpha ud-\phi)=0$, which implies that $\A$ is not a domain and the elements $ud, du, t_{1}, \cdots, t_{n}$ are not algebraically independent. Furthermore, we can show that the algebra $\A$ is neither a right nor a left noetherian using the same construction of the infinite ascending chain of ideals as in {\bf Lemma 4.3} in \cite{KMP}. Thus we have proved the implications: $(1)\Rightarrow (3), (2)\Rightarrow (3)$ and $(4)\Rightarrow (3)$. 

When $\beta\neq 0$, we can show that the graded algebra $gr(\A)$ of $\A$ is indeed Auslander-regular. As a result, the algebra $\A$ itself is Auslander-regular. \qed

In addition, we have the following results which can be directly derived from the proof of the previous theorem.
\begin{prop}
Assume that $\beta \neq 0$ and let $R=\K[x, y, t_{1}, \cdots, t_{n}]$. Then the noetherian down-up algebra $\A$ is isomorphic to a generalized Weyl algebra $R(\sigma,x)$.
\end{prop}
\qed

\begin{prop}
Assume that $\beta \neq 0$ and $\K$ is sufficiently large, then $\A$ has a filtration such that the graded algebra $gr(\A)$ is isomorphic to an iterated skew polynomial ring of the form:
\[
\K[t_{1}, \cdots, t_{n}][x_{1}][x_{2};\eta][x_{3};\tau, \delta].
\]
\end{prop}\qed

Moreover, we also have the following proposition.
\begin{prop}
When $\beta\neq 0$, any invertible element in $\A$ belongs to $\K^{\ast}$.
\end{prop}
\qed

Following the proof of {\bf Lemma 1.2} in \cite{CM}, we can similarly prove the following result.
\begin{prop}
When $\beta \neq 0$, the set $\mathbb{D}=\{d^{i}\mid i\geq 0\}$ (resp. $\mathbb{U}=\{u^{i}\mid i\geq 0\}$ ) is a (left and right) Ore set of $\A$. And the localization $\A_{\mathbb{D}}$ (resp. $\A_{\mathbb{U}}$) is isomorphic to the skew group ring $R[z^{\pm 1}, \sigma]$ (resp. $R[z^{\pm 1}, \sigma^{\prime}]$ where $\sigma^{\prime}$ is similarly defined as $\sigma$) with $R=\K[ud, du, t_{1}, \cdots, t_{n}]$. 
\end{prop}

Additionally, we have the following result about the left global dimension of $\A$.
\begin{thm}
The down-up algebra $\A$ has a global dimension of $n+3$.
\end{thm}
{\bf Proof:} First of all, the down-up algebra $\A$ has a filtration via assigning $deg(u)=deg (d)=1$ and $deg(t_{1})=\cdots=deg(t_{n})=0$ such that the corresponding graded algebra ${\rm gr}(\A)$ is isomorphic to an iterated skew polynomial ring $R[t_{1}, \cdots, t_{n}]$ where $R=A(\alpha, \beta, 0)$ is a ``classical" graded down-up algebra. Note that we have ${\rm gl dim}(R)=3$ from \cite{KK2, KMP}. Therefore, by {\bf Theorem 7.5.3 \& Corollary 7.6.18} in \cite{MR}, we have that ${\rm gl dim} (\A)\leq {\rm gldim}( {\rm gr}(\A))={\rm gldim}(R)+n=n+3$. 

Second of all, we also have that $\A/(t_{1}-\lambda_{1})=\mathcal{A}(\alpha, \beta, \phi^{\prime})$ for some $\phi^{\prime} \in \K[t_{2}, \cdots, t_{n}]$, which is thus a down-up algebra defined over $\K[t_{2}, \cdots, t_{n}]$. Since the element $t_{1}-\lambda_{1}$ is regular in the algebra $\A$, by {\bf Theorem 7.3.5} in \cite{MR}, we have that ${\rm gldim}(\A)\geq {\rm gldim}(\mathcal{A}(\alpha, \beta, \phi^{\prime}))+1$. Therefore, using induction, we can now show that ${\rm gldim}(\A)\geq n+3$. So we have proved that the global dimension of $\A$ is $n+3$. \qed

It is easy to see that the down-up algebra $\A$ is rarely a primitive ring. However, the algebra $\A$ is still a prime ring in most of the cases. 
\begin{thm}
The algebra $\A$ is a prime ring except in the case when $\alpha=\beta=\phi=0$.
\end{thm}
{\bf Proof:} First of all, we know that the multiplicative set 
\[
\mathbb{S}=\{f(t_{1}, \cdots, t_{n})\mid 0\neq f(t_{1}, \cdots, t_{n})\in \K[t_{1}, \cdots, t_{n}]\}
\]
consists of regular central elements of $\A$. Thus the down-up algebra $\A$ is a prime ring if and only if its localization $\A_{\mathbb{S}}$ is a prime ring. Now the result follows from {\bf Theorem 3.2} in \cite{KK2}. \qed

Recall that the Krull dimensions of the ``classical" noetherian down-up algebras $A(\alpha, \beta, \gamma)$ were studied in \cite{BL, BV}. In particular, the Krull dimensions of the ``classical" noetherian down-up algebras were explicitly calculated in \cite{BL}. Based on the results in \cite{BL}, we shall have the following result on the Krull dimension of the noetherian down-up algebra $\A$ (i.e., we will assume that $\beta \neq 0$). 
\begin{thm} We have the following description on the Krull dimension of $\A$:
\begin{enumerate}
\item If $\phi$ is a constant, then $\A$ has a Krull dimension of $n+2$ if and only if $char \K=0, \gamma\neq 0$, and $\alpha+\beta=1$; otherwise, the Krull dimension of $\A$ is $n+3$;\\

\item If $\phi$ is not a constant and $\K$ is an algebraically closed field, then the Krull dimension of $\A$ is $ n+3$.
\end{enumerate}
\end{thm}
{\bf Proof:} If $\phi$ is a constant, then the algebra $\A$ is isomorphic to the tensor product $A(\alpha, \beta, \gamma)\otimes_{\K} \K[t_{1}, \cdots, t_{n}]$ of a ``classical" down-up algebra $A(\alpha, \beta, \gamma)$ where $\gamma=\phi$ with the polynomial ring $\K[t_{1}, \cdots, t_{n}]$ over the base field $\K$. Note that {\bf Theorem 4.2} in \cite{BL} states that the ``classical" down-up algebra $A(\alpha, \beta, \gamma)$ has a Krull dimension of $2$ if and only if $char\K=0, \gamma \neq 0$, and $\alpha+\beta=1$; otherwise, the Krull dimension of $A(\alpha, \beta, \gamma)$ is $3$. Thus, by {\bf Proposition 6.5.4 } in \cite{MR}, we can prove that the noetherian down-up algebra $\A$ has a Krull dimension of $n+2$ if and only if $char \K=0, \phi \neq 0$, and $\alpha+\beta=1$; otherwise, it has a Krull dimension of $n+3$. 

If $\phi$ is not a constant and $\K$ is an algebraically closed field, then we can choose scalars $\lambda_{1}, \cdots, \lambda_{n}\in \K$ such that $\phi(\lambda_{1}, \cdots, \lambda_{n})=0$. It is obvious that $\A/(t_{1}-\lambda_{1}, \cdots, t_{n}-\lambda_{n})$ is isomorphic to a ``classical" 
down-up algebra $A(\alpha, \beta, 0)$ defined over the field $\K$. Thus, via {\bf Proposition 6.3.11} in \cite{MR} and using induction, we have that $\mathcal{K}(\A)\geq n+3$. Recall that $\A$ has a filtration such that the corresponding graded algebra ${\rm gr}(\A)$ is isomorphic to $A(\alpha, \beta, 0)\otimes_{\K}\K[t_{1}, \cdots, t_{n}]$. Using {\bf Proposition 6.5.4 \& Lemma 6.5.6} in \cite{MR}, we have that $\mathcal{K}(\A)\leq \mathcal{K}(A(\alpha, \beta, 0))+n=3+n$. So we have proved that $\mathcal{K}(\A)=n+3$ when $\phi$ is not a constant and $\K$ is an algebraically closed field.\qed

\section{Automorphisms and Isomorphisms of $\A$}

When a family of algebras is defined by generators and relations involving parameters, one is led to study the classification of these algebras up to isomorphisms for the purpose of understanding the dependence on the parameters. In addition, it is also important to study the symmetries of an algebra in terms of its algebra automorphisms. Note that the isomorphism problems and automorphism groups of some classes of generalized Weyl algebras were studied in \cite{BJ} and the references therein. The isomorphism classification of down-up algebras was initiated in \cite{BR} via studying one-dimensional modules of $\A$, and the isomorphism problem was nearly perfectly settled in \cite{CM}. The automorphism group of a ``classical" down-up algebra $A(\alpha, \beta, \gamma)$ was determined \cite{CL} under the assumption that either $r$ or $s$ is not a root of unity.

In this section, we will study the automorphisms and isomorphisms of the noetherian down-up algebra $\A$ under some restrictions on the parameters. Due to the complexity of the problem, we are only able to establish some results on isomorphisms and automorphisms of $\A$ for certain special cases. From now on, we will assume that the field $\K$ is of characteristic zero and large enough so that there exist $r,s\in \K^{\ast}$, which are not roots of unity such that $r+s=\alpha$ and $rs=-\beta$. The approach is to first determine the normal elements of the algebra $\A$ and study the images of normal elements under any algebra automorphism of $\A$. Though some of the results could be somehow strengthened to cover more cases, we will focus on two major cases where $r,s$ are either not related in the sense that $r^{i}s^{j}=1$ implies $i=j=0$, or particularly $r=s^{-1}$. 

First of all, we are going to rewrite the last two defining relations of the algebra $\A$ using the new parameters $r,s$ as follows:
\begin{eqnarray*}
d^{2}u-(r+s)dud+rsud^{2}-\phi d=0;\\
du^{2}-(r+s)udu+rsu^{2}d-\phi u=0.
\end{eqnarray*}

Let us set $H=du-rud+\frac{\phi}{s-1}$ and $K=du-sud+\frac{\phi}{r-1}$. Then we have the following:
\[
dH=sHd, \quad Hu=suH;\quad dK=rKd,\quad Ku=ruK.
\]
It is also straightforward to verify that the set 
\[
\{H^{i}K^{j}u^{k}t_{1}^{m_{1}} \cdots t_{n}^{m_{n}}, H^{i}K^{j}d^{k+1}t_{1}^{m_{1}}\cdots t_{n}^{m_{n}}\mid i, j, k, m_{1}, \cdots, m_{n}\geq 0\}
\]
is also a $\K-$basis of the algebra $\A$. As a matter of fact, we can also consider the algebra $\A$ as a $\mathbb{Z}-$graded algebra with 
\begin{eqnarray*}
\A_{0}=\K[H, K, t_{1}, \cdots, t_{n}];\\
\A_{i}=\K [H, K, t_{1}, \cdots, t_{n}] d^{i},\\
\A_{-i}=\K[H, K, t_{1}, \cdots, t_{n}] u^{i}
\end{eqnarray*}
for any $i\geq 1$.
\begin{lem}
Let $a\in \A$ be a normal element of $\A$ and $0\neq f(t_{1}, \cdots, t_{n}) \in \K[t_{1}, \cdots, t_{n}]$. Then the element $a/f(t_{1}, \cdots, t_{n})$ is a normal element of the localization $\A_{\mathbb{S}}$. Conversely, any normal element of $\A_{\mathbb{S}}$ is of the form $a/f(t_{1}, \cdots, t_{n})$ for some normal element of $\A$ and a non-zero polynomial $f(t_{1}, \cdots, t_{n})\in \K[t_{1}, \cdots, t_{n}]$.
\end{lem}
{\bf Proof:} Note that any non-zero element $f(t_{1}, \cdots, t_{n})\in \K[t_{1}, \cdots, t_{n}]$ is a regular central element of $\A$ and thus the multiplicative set $\mathbb{S}=\{f(t_{1}, \cdots, t_{n})\mid 0\neq f(t_{1}, \cdots, t_{n})\in \K[t_{1}, \cdots, t_{n}]\}$ consists of non-zero central regular elements of $\A$. Therefore, $\A$ can be regarded as a subalgebra of $\A_{\mathbb{S}}$; and the result follows.\qed

\begin{prop}
\begin{enumerate}
\item If $r^{i}s^{j}=1$ implies that $i=j=0$, then any normal element of the algebra $\A$ is of the form $f(t_{1}, \cdots, t_{n})H^{i}K^{j}$ for some $f(t_{1}, \cdots, t_{n})\in \K[t_{1}, \cdots, t_{n}]$ and $i, j\geq 0$.\\

\item If $r=s^{-1}$, then any normal element of the algebra $\A$ is of the form $g_{1}(HK, t_{1}, \cdots, t_{n})H^{i}$ or $g_{2}(HK, t_{1}, \cdots, t_{n})K^{j}$ where $g_{1}(HK, t_{1}, \cdots, t_{n})$, $g_{2}(HK, t_{1}, \cdots, t_{n})\in \K[HK, t_{1}, \cdots, t_{n}]$ and $i, j\geq 0$. 
\end{enumerate}
\end{prop}
{\bf Proof:} First of all, we prove that the corresponding result holds for the ``classical" down-up algebra $\A_{\mathbb{S}}$, which is the localization of $\A$ at the multiplicative set $\mathbb{S}$. Suppose that $N\in \A_{\mathbb{S}}$ is a normal element. Using {\bf Proposition 2.5} in \cite{CL}, we know that each normal element of the algebra $\A_{\mathbb{S}}$ is of the form $f(H, K)x^{n}$ or $f(H, K)y^{n}$ for some $f(H, K)\in \K(t_{1}, \cdots, t_{n})[H, K]$ and $n\geq 0$ such that $f(H,K), x^{n}$ and $y^{n}$ are normal elements of $\A_{\mathbb{S}}$. Since $r^{i}s^{j}=1$ implies that $i=j=0$, we know that the elements $x^{n}$ and $y^{n}$ with $n\geq 1$ cannot be normal elements of the algebra $\A_{\mathbb{S}}$ using {\bf Lemma 2.4} in \cite{CL}. Therefore, we have that $N\in (\A_{\mathbb{S}})_{0}=\K(t_{1}, \cdots, t_{n})[H, K]$. Suppose that $N=\sum_{i, j} f_{ij}(t_{1}, \cdots, t_{n}) H^{i}K^{j}$. Then we have the following:
\begin{eqnarray*}
(\sum_{ij} f_{ij}(t_{1},\cdots, t_{n}) H^{i}K^{j})d=d(\sum_{ij} f_{ij}(t_{1}, \cdots, t_{n}) r^{i}s^{j} H^{i}K^{j}),\\
(\sum_{ij} f_{ij}(t_{1}, \cdots, t_{n}) H^{i}K^{j})u=u(\sum_{ij} f_{ij}(t_{1},\cdots, t_{n})r^{-i}s^{-j} H^{i}K^{j})
\end{eqnarray*}
which implies that $r^{i}s^{j}=\lambda$ and $r^{-i}s^{-j}=\mu$ for all $i, j$. Since $r^{i}s^{j}=1$ implies $i=j=0$, we have $f_{ij}(t_{1},\cdots, t_{n})=0$ except for one pair of $i, j$. So we have proved that $N=f_{ij}(t_{1}, \cdots, t_{n})H^{i}K^{j}$ for some $i, j$. Via {\bf Lemma 3.1}, we have the desired result for the algebra $\A$. Thus, we have proved the statement in $(1)$. 

The statement in $(2)$ can be proved in a similar fashion and we will not repeat the details here. \qed

Now we assume that $n=1$ and denote $t_{1}$ by $t$. We have the following description of the automorphisms of the algebra $\A$ under some conditions on $r,s$.
\begin{thm}
\begin{enumerate}
\item If $r^{i}s^{j}=1$ implies that $i=j=0$, then any algebra automorphism $\sigma\colon \A\longrightarrow \A$ is of the form: 
\[
\sigma(d)=\lambda_{1}d,\quad \sigma(u)=\lambda_{2} u,\quad \sigma(t)=at+b
\]
for some $\lambda_{1}, \lambda_{2}, a\in \K^{\ast}, b\in \K$ such that $\lambda_{1}\lambda_{2} \phi(t)=\phi (at+b)$. In particular, if $\phi=0$, then we have the following:
\[
\sigma(d)=\lambda_{1} d, \quad \sigma(u)=\lambda_{2} u, \quad \sigma(t)=at+b
\]
for some $\lambda_{1}, \lambda_{2}, a\in \K^{\ast}$ and $b\in K$; if $\phi$ is a non-zero constant, then we have the following:
\[
\sigma(d)=\lambda d,\quad \sigma(u)=\lambda^{-1} u,\quad \sigma(t)=at+b
\]
for some $\lambda, a\in \K^{\ast}$ and $b\in \K$; if $\phi=\sum_{i=0}^{m}a_{i}t^{i}$ is not a constant, then we have the following:
\[
\sigma(d)=\lambda_{1}, \quad \sigma(u)=\lambda_{2} u,\quad \sigma(t)=at
\]
for some $\lambda_{1}, \lambda_{2}, a \in \K^{\ast}$ such that $a^{i}=\lambda_{1}\lambda_{2}$ whenever $a_{i}\neq 0$.\\

\item If $r=s^{-1}$, then any algebra automorphism $\sigma$ of $\A$ is of the form:
\[
\sigma(d)=\lambda_{1} d, \quad \sigma(u)=\lambda_{2}u,\quad \sigma(t)=at+g(HK)
\]
or
\[
\sigma(d)=\lambda_{1} u, \quad \sigma(u)=\lambda_{2}d,\quad \sigma(t)=at+g(HK)
\]
for some $\lambda_{1}, \lambda_{2}, a\in \K^{\ast}$ and $g(HK)\in \K[HK]$ such that $\lambda_{1}\lambda_{2}\phi(t)=\phi(at+g(HK))$. In particular, if $\phi=0$, then we have the following:
\[
\sigma(d)=\lambda_{1} d, \quad \sigma(u)=\lambda_{2}u,\quad \sigma(t)=at+g(HK)
\]
or 
\[
\sigma(d)=\lambda_{1} u, \quad \sigma(u)=\lambda_{2} d, \quad \sigma(t)=at+g(HK)
\]
for some $\lambda_{1}, \lambda_{2}, a \in \K^{\ast}$ and $g(HK)\in \K[HK]$; if $\phi$ is a non-zero constant, then we have the following:
\[
\sigma(d)=\lambda u, \quad \sigma(u)=\lambda^{-1}u,\quad \sigma(t)=at+g(HK)
\]
or 
\[
\sigma(d)=\lambda u, \quad \sigma(u)=\lambda^{-1}d, \quad \sigma(t)=at+g(HK)
\]
for some $\lambda, a\in \K^{\ast}$ and $g(HK)\in \K[HK]$; if $\phi=\sum_{i=1}^{m} a_{i}t^{i}$ is not a constant, then we have the following:
\[
\sigma(d)=\lambda_{1} d, \quad \sigma(u)=\lambda_{2} u,\quad \sigma(t)=at
\]
or
\[
\sigma(d)=\lambda_{1}u, \quad \sigma(u)=\lambda_{2} d,\quad \sigma(t)=at
\]
for some $\lambda_{1}, \lambda_{2}, a \in \K^{\ast}$ such that $a^{i}=\lambda_{1}\lambda_{2}$ whenever $a_{i}\neq 0$.\\
\end{enumerate}
\end{thm}
{\bf Proof:} Let $\sigma \colon \A\longrightarrow \A$ be an algebra automorphism of $\A$. Since $H$ and $K$ are normal elements of $\A$, the images of $H$ and $K$ under $\sigma$ are also normal elements of $\A$. Therefore, we have $\sigma(H)=f_{1}H^{a}K^{b}$ and $\sigma(K)=f_{2}H^{c}K^{d}$ for some $f_{1}, f_{2}\in \K[t]$ and $a, b, c, d \in \mathbb{Z}_{\geq 0}$. In addition, we also have $\sigma^{-1}(H)=g_{1}H^{x}K^{y}$ and $\sigma^{-1}(K)=g_{2}H^{z}K^{w}$ for some $g_{1}, g_{2}\in \K[t]$ and $x, y, z, w\in \mathbb{Z}_{\geq 0}$. Applying the inverse of $\sigma$ to the elements $\sigma(H)$ and $\sigma(K)$, we shall have the following:
\[
H=\sigma^{-1}(\sigma(H))=\sigma^{-1}(f_{1}H^{a}K^{b})=\sigma^{-1}(f_{1}) (g_{1}H^{x}K^{y})^{a}(g_{2}H^{z}K^{w})^{b}
\]
and 
\[
K=\sigma^{-1}(\sigma(f_{2}H^{c}K^{d}))=\sigma^{-1}(f_{2}H^{c}K^{d})=\sigma^{-1}(f_{2}) (g_{1}H^{x}K^{y})^{c}(g_{2}H^{z}K^{w})^{d}.
\]
As a result, we have the following:
\[
\sigma^{-1}(f_{1})g_{1}^{a}g_{2}^{b}=1,\quad \sigma^{-1}(f_{2})g_{1}^{c}g_{2}^{d}=1;
\]
and
\[
H^{ax+bz}K^{ay+bw}=H,\quad H^{cx+dz}K^{cy+dw}=K.
\]
Therefore, we have $f_{1}, f_{2}, g_{1}, g_{2}\in K^{\ast}$ and either $a=1, b=0, c=0, d=1$ or $a=0, b=1, c=1, d=0$. Since $r^{i}s^{j}=1$ implies $i=j=0$, we can further show that $a=1, b=0, c=0, d=1$ using the commuting relations between $u,d$ and $H, K$. Thus we have $\sigma(H)=\mu_{1} H$ and $\sigma (K)=\mu_{2}K$. As a result, we can also show that $\sigma(d)=\lambda_{1} d$ and $\sigma(u)=\lambda_{2} u$. Since the center of $\A$ is just $\K[t]$ and the center is stable under any algebra automorphism, we have $\phi(t)=at+b$ for some $a\in \K^{\ast}, b\in \K$. Now we need to determine the relationship among the scalars $\lambda_{1}, \lambda_{2}, a, b$. Applying $\sigma$ to the defining relations, we shall have the following:
\[
\lambda_{1}^{2}\lambda_{2}d^{2}u-\lambda_{1}^{2}\lambda_{2}(r+s)dud+\lambda_{1}^{2}\lambda_{2}rsud^{2}-\lambda_{1} \phi(at+b)d=0,
\]
and 
\[
\lambda_{1}\lambda_{2}^{2}du^{2}-\lambda_{1}\lambda_{2}^{2}(r+s)udu+\lambda_{1}\lambda_{2}^{2} rsu^{2}d-\lambda_{2}\phi(at+b)u=0.
\]

Thus we have $\lambda_{1}\lambda_{2}\phi(t)=\phi (at+b)$. And any automorphism $\sigma$ of $\A$ is defined by 
\[
\sigma(d)=\lambda_{1}d,\quad \sigma(u)=\lambda_{2} u,\quad \sigma(t)=at+b
\]
for some $\lambda_{1}, \lambda_{2}, a\in \K^{\ast}, b\in \K$ such that $\lambda_{1}\lambda_{2} \phi(t)=\phi (at+b)$. Conversely, any mapping $\sigma$ defined as above can be extended to an algebra automorphism of $\A$. In addition, if $\phi(t)=0$, then $\sigma(d)=\lambda_{1} d, \sigma(u)=\lambda_{2} u, \sigma(t)=at+b$ for any $\lambda_{1}, \lambda_{2}, a \in \K^{\ast}$ and $b\in \K$. If $\phi$ is a non-zero constant, then $\sigma(d)=\lambda d, \sigma(u)=\lambda^{-1}u, \sigma(t)=at+b$ for some $\lambda, a \in \K^{\ast}$ and $b\in \K$. If $\phi=\sum_{i=1}^{m} a_{i}t^{i}$ is not a constant, then we have the following:
\[
\sigma(d)=\lambda_{1} d, \quad \sigma(u)=\lambda_{2} u,\quad \sigma(t)=at
\]
or
\[
\sigma(d)=\lambda_{1}u, \quad \sigma(u)=\lambda_{2} d,\quad \sigma(t)=at
\]
for some $\lambda_{1}, \lambda_{2}, a \in \K^{\ast}$ such that $a^{i}=\lambda_{1}\lambda_{2}$ whenever $a_{i}\neq 0$. Thus we have proved $(1)$.

Now we assume that $r=s^{-1}$. Let $\sigma\colon \A \longrightarrow \A$ be an algebra automorphism of $\A$. Similarly, one can show that $\sigma(H)=\mu_{1}H, \sigma(K)=\mu_{2}K$ or $\sigma(H)=\mu_{1}K, \sigma(K)=\mu_{2}H$ for some $\mu_{1}, \mu_{2}\in \K^{\ast}$. Using the commuting relationship between $d, u$ and $H,K$, one can show that if $\sigma(H)=\mu_{1}H, \sigma(K)=\mu_{2}K$, then $\sigma(d)=\lambda_{1}d, \sigma(u)=\lambda_{2} u$ for some $\lambda_{1}, \lambda_{2}\in \K^{\ast}$; and if $\sigma(H)=\mu_{1}K, \sigma(K)=\mu_{2}H$, then $\sigma(d)=\lambda_{1}u, \sigma(u)=\lambda_{2}d$ for some $\lambda_{1}, \lambda_{2}\in \K^{\ast}$.
In either case, we shall have $\sigma(t)=at+g(HK)$ for some $a\in K^{\ast}, g(HK)\in \K[HK]$ because the center of $\A$ is $\K[HK, t]$ and $\sigma(HK)=\mu_{1}\mu_{2}HK$. If $\sigma(H)=\mu_{1}H, \sigma(K)=\mu_{2}K$, then we have the following:
\[
\sigma(d)=\lambda_{1} d, \quad \sigma(u)=\lambda_{2}u,\quad \sigma(t)=at+g(HK)
\]
for some $\lambda_{1}, \lambda_{2}, a\in \K^{\ast}$ and $g(HK)\in \K[HK]$ such that $\lambda_{1}\lambda_{2}\phi(t)=\phi(at+g(HK))$. If $\sigma(H)=\mu_{1}K, \sigma(K)=\mu_{2}H$, then we have the following:
\[
\sigma(d)=\lambda_{1} u, \quad \sigma(u)=\lambda_{2}d,\quad \sigma(t)=at+g(HK)
\]
for some $\lambda_{1}, \lambda_{2}, a\in \K^{\ast}$ and $g(HK)\in \K[HK]$ such that $\lambda_{1}\lambda_{2}\phi(t)=\phi(at+g(HK))$ . Therefore, when $r=s^{-1}$, any algebra automorphism $\sigma$ of $\A$ is of the form:
\[
\sigma(d)=\lambda_{1} d, \quad \sigma(u)=\lambda_{2}u,\quad \sigma(t)=at+g(HK)
\]
or
\[
\sigma(d)=\lambda_{1} u, \quad \sigma(u)=\lambda_{2}d,\quad \sigma(t)=at+g(HK)
\]
for some $\lambda_{1}, \lambda_{2}, a\in \K^{\ast}$ and $g(HK)\in \K[HK]$ such that $\lambda_{1}\lambda_{2}\phi(t)=\phi(at+g(HK))$. Conversely, any mapping $\sigma$ as defined above can be extended to an algebra automorphism of $\A$. Additionally, one can similarly verify the result on $\lambda_{1}, \lambda_{2}$, and $g(HK)$ when $\phi$ takes special values accordingly. Thus we have proved $(2)$. \qed

As an application of the proof of the previous theorem, we shall also have the following result on the isomorphisms between down-up algebras $\A$.

\begin{thm}
\begin{enumerate}
\item If $r_{1}^{i}s_{2}^{j}=1$ implies $i=j=0$ and $r_{2}^{i}s_{2}^{j}=1$ does not imply $i=j=0$, then $\mathcal{A}(r_{1},s_{1}, \phi_{1})$ is never isomorphic to $\mathcal{A}(r_{2},s_{2}, \phi_{2})$.\\

\item If $r^{i}s^{j}=1$ implies $i=j=0$, then $\mathcal{A}(r,s, \phi_{1}) \cong \mathcal{A}(r,s, \phi_{2})$ if and only if $\phi_{1}(at+b)=\eta \phi_{2}(t)$ for some $\eta \in \K^{\ast}$.\\

\item If $r_{1}^{i_{1}}s_{1}^{j_{1}}=r_{2}^{i_{2}}s_{2}^{j_{2}}=1$ implies $i_{1}=j_{1}=i_{2}=j_{2}=0$, then $\mathcal{A}(r_{1},s_{1}, \phi_{1})$ is isomorphic to $\mathcal{A}(r_{2}, s_{2}, \phi_{2})$ if and only if\\

\begin{enumerate}
\item either $r_{1}=r_{2}, s_{1}=s_{2}$ and $\phi_{1}(at+b)=\eta \phi_{2}(t)$ for some $\eta, a\in \K^{\ast}$ and $b\in \K$;\\

\item or $r_{1}=s_{2}, r_{2}=s_{1}$ and $\phi_{1}(at_{2}+b)=\eta \phi_{2}(t)$ for some $\eta, a\in \K^{\ast}$ and $b\in \K$;\\

\item or $r_{1}=s_{2}^{-1}, r_{2}=s_{1}^{-1}$ and $\phi_{1}(at+b)=\eta \phi_{2}(t)$ for some $\eta, a\in \K^{\ast}$ and $b\in \K$;\\

\item or $r_{1}=r_{2}^{-1}, s_{1}=s_{2}^{-1}$ and $\phi_{1}(at+b)=\eta \phi_{2}(t)$ for some $\eta, a\in \K^{\ast}$ and $b\in \K$.\\

\end{enumerate}

\item $\mathcal{A}(r,r^{-1}, \phi_{1})\cong \mathcal{A}(r,r^{-1}, \phi_{2})$ if and only if $\phi_{1}(at+b)=\eta \phi_{2}(t)$ for some $\eta, a\in \K^{\ast}$ and $b \in \K$. 
\end{enumerate}
\end{thm}
{\bf Proof:} Part $(1)$ follows from the fact that the algebra $\mathcal{A}(r_{1}, s_{1}, \phi_{1})$ has a center $\K[t]$, while the algebra $\mathcal{A}(r_{2}, s_{2}, \phi_{2})$ has a center strictly containing $\K[t]$. 

When $r^{i}s^{j}=1$ implies $i=j=0$, both the algebra $\mathcal{A}(r, s, \phi_{1})$ and the algebra $\mathcal{A}(r, s, \phi_{2})$ have a center equal to $\K[t]$ respectively. Thus, the set of the normal elements of $\mathcal{A}(r,s, \phi_{1})$ is given by the set $\{g(t)H^{i}K^{j}\mid i, j\geq 0, g(t) \in \K[t]\}$ and the set of the normal elements of $\mathcal{A}(r,s, \phi_{2})$ is given by the set $\{g(t) H^{i}K^{j}\mid i, j\geq 0, g(t) \in \K[t]\}$. If $\sigma \colon \mathcal{A}(r,s, \phi_{1})\longrightarrow \mathcal{A}(r,s, \phi_{2})$ is an algebra isomorphism, then $\sigma(H)=\mu_{1} H$ and $\sigma(K)=\mu_{2}K$ and $\sigma(t)=at+b$ for some $\mu_{1}, \mu_{2}, a\in K^{\ast}, b\in \K$. As a result, we have $\sigma(d)=\lambda_{1}d$, $\sigma(u)=\lambda_{2}u$ and $\sigma(t)=at+b$ for some $\lambda_{1}, \lambda_{2}, a\in \K^{\ast}, b\in \K$. Applying the isomorphism $\sigma$ to the defining relations, we shall have $\lambda_{1}\lambda_{2}\phi_{2}(t)=\phi_{1}(at+b)$. Thus we have $\eta \phi_{2}(t)=\phi_{1}(at+b)$ for some $\eta, a\in \K^{\ast}, b\in \K$. Conversely, if $\eta \phi_{2}(t)=\phi_{1}(at+b)$ for some $\eta, a\in \K^{\ast}, b\in \K$, then it is obvious that one can define an algebra isomorphism $\sigma$ from $\mathcal{A}(r, s, \phi_{1})$ onto $\mathcal{A}(r, s, \phi_{2})$. Thus we have proved part $(2)$.

Similarly, one can prove the statements in parts (3) and (4), and we will not state the details here. \qed

\end{document}